\documentclass[12pt,a4paper]{article}

\usepackage[T1]{fontenc}
\usepackage[utf8]{inputenc}
\usepackage{lmodern}
\usepackage{amsmath,amssymb,amsthm}
\usepackage{geometry}
\usepackage{microtype}
\usepackage{array}
\usepackage{float}
\usepackage{enumitem}
\usepackage{cite}
\usepackage{hyperref}
\usepackage{orcidlink}

\hypersetup{
  colorlinks=true,
  linkcolor=blue,
  citecolor=blue,
  urlcolor=blue,
  pdftitle={Universal Invariants for gcd(n\string^k-1,n!-1)>1},
  pdfauthor={Marco Rip\`a},
  pdfsubject={Factorial residues and greatest common divisors},
  pdfkeywords={factorial residues, greatest common divisors, multiplicative orders, roots of unity, universal invariants}
}

\newtheorem{definition}{Definition}[section]
\newtheorem{lemma}[definition]{Lemma}
\newtheorem{proposition}[definition]{Proposition}
\newtheorem{theorem}[definition]{Theorem}
\newtheorem{corollary}[definition]{Corollary}
\newtheorem{conjecture}[definition]{Conjecture}
\newtheorem{problem}[definition]{Open Problem}
\theoremstyle{remark}
\newtheorem{remark}[definition]{Remark}

\newcommand{\N}{\mathbb{N}}
\newcommand{\F}{\mathbb{F}}
\newcommand{\inv}{\operatorname{inv}}
\newcommand{\oddomega}{\omega_{\mathrm{odd}}}

\newcommand{\keywords}[1]{%
  \par
  {\leftskip=2.5em \rightskip=2.5em
  \noindent\textbf{Keywords:} #1\par}
  \smallskip
}

\newcommand{\subjclass}[1]{%
  \par
  {\leftskip=2.5em \rightskip=2.5em
  \noindent\textbf{2020 Mathematics Subject Classification:} #1\par}
}

\newcommand{\doi}[1]{\href{https://doi.org/#1}{\nolinkurl{doi:#1}}}

\begin{document}

\title{Universal Invariants for $\gcd(n^k-1,n!-1)>1$
\vspace{1.5em}}
\author{Marco Rip\`a~\orcidlink{0000-0002-6036-5541}\\[0.4ex]\normalsize Rome, Italy}
\date{}

\maketitle
\vspace{1em}

\begin{abstract}
For integers $k\geq2$ and $\lambda_1,\lambda_2\in\{-1,+1\}$, let $S_k^{\lambda_1,\lambda_2}:=\{n \in \mathbb{N}\setminus \{0,1\} :\gcd(n^k+\lambda_1,n!+\lambda_2)>1\}$.
This family arose from a MathOverflow problem concerning the existence of integers $n$ satisfying $\gcd(n^2+1,n!+1)>1$ and from the subsequent question of whether arbitrarily many, but finitely many, solutions can be obtained by varying the exponent and the signs.
We concentrate on the case $\lambda_1=\lambda_2=-1$ and introduce the notion of a universal invariant for an odd prime $p$, namely an integer $n\geq2$ such that $\gcd(n^{p \cdot m}-1,n!-1)>1$ for every integer $m\geq1$.
We prove that this condition is equivalent to $\gcd(n^p-1,n!-1)>1$, characterize the possible prime divisors in terms of roots of unity over finite fields, and recall an explicit construction due to Jonathan Love that yields a universal invariant for every odd prime $p\geq5$, with the cubic case treated separately. As an application, we prove that $|S_k^{-1,-1}|\geq\oddomega(k)$ for every integer $k\geq2$, where $\oddomega(k)$ denotes the number of distinct odd prime divisors of $k$.
Hence, the cardinalities of these sets are unbounded as $k$ varies.
This lower bound is deliberately coarse and constructive, as it is obtained by considering only universal invariants; stronger general lower bounds are obtained by different constructions.
We also discuss minimal universal invariants, their persistent common divisors, the special structure of the cubic case, and a heuristic due to Will Sawin suggesting that there may be infinitely many cubic universal invariants, although the known examples are exceptionally sparse. Whether $S_k^{-1,-1}$ is finite or infinite for any fixed exponent $k\geq3$ remains an open problem.
\end{abstract}

\keywords{Factorial residues, Greatest common divisors, Multiplicative orders, Roots of unity, Universal invariants.}

\subjclass{Primary 11A07; Secondary 11A41, 11B65.}
\vspace{1.5em}

\section{Introduction}
We study common divisors of powers and factorials through the family of inequalities
\[
\gcd(n^k+\lambda_1,n!+\lambda_2)>1,
\qquad
n,k\geq2,
\qquad
\lambda_1,\lambda_2\in\{-1,+1\},
\]
where $n$ and $k$ are both integers.

Our main focus is the negative-sign family $\gcd(n^k-1,n!-1)>1$, for which persistent solutions associated with odd prime divisors of the exponent lead to the universal invariants studied in the present paper and to constructive lower bounds for the number of solutions.

The broader problem originated from a MathOverflow question posted by the user ``smaller'' \cite{SmallerMO}, asking whether
\[
\gcd(n^2+1,n!+1)=1
\]
holds for every integer $n\geq2$. Computations up to $4\cdot10^6$ had produced no counterexample (see comments in \cite{SmallerMO}), and the answers established nontrivial restrictions on any possible solution. In particular, Adrien Mounier showed that a nontrivial gcd would necessarily be prime \cite{MounierMO}. Nevertheless, the existence of an integer $n$ satisfying
\[
\gcd(n^2+1,n!+1)>1
\]
remains open.

This question naturally leads to varying both the exponent and the signs. Accordingly, for fixed $k\geq2$ and $\lambda_1,\lambda_2\in\{-1,+1\}$, define
\[
S_k^{\lambda_1,\lambda_2}
:=
\{n\in\N\setminus\{0,1\}:
\gcd(n^k+\lambda_1,n!+\lambda_2)>1\}.
\]
This broader line of investigation was developed by the author through a
sequence of MathOverflow discussions
\cite{RipaQuadraticMO,RipaCardinalityMO,RipaCubicMO,RipaUniversalMO}. Its original motivation was not merely to decide whether an isolated set $S_k^{\lambda_1,\lambda_2}$ was empty or nonempty, but to understand the possible cardinalities of these sets. More precisely, the question posed in \cite{RipaCardinalityMO} was whether, for every positive integer $c$, there exists a triple $(k,\lambda_1,\lambda_2)$ such that
\[
c<|S_k^{\lambda_1,\lambda_2}|<+\infty.
\]
Thus, the intended phenomenon was intermediate between isolated solvability and an immediately infinite family.

The first cases revealed very different behaviours. Marco Mantovanelli
proved that
\[
S_2^{-1,-1}=\varnothing,
\]
or equivalently
\[
\gcd(n^2-1,n!-1)=1
\]
for every $n\geq2$ \cite{MantovanelliMO}. At the opposite extreme, $S_3^{+1,+1}$ is infinite. Indeed, if $q>2$ is prime and $n=q-1$, then
\[
n^3+1\equiv(-1)^3+1\equiv0\pmod q,
\]
while Wilson's theorem gives
\[
n!+1=(q-1)!+1\equiv0\pmod q.
\]
Consequently, every odd prime occurs as a common divisor, and every integer one less than an odd prime belongs to $S_3^{+1,+1}$. The investigation then turned to the negative-sign family
\[
S_k^{-1,-1}
=
\{n \in \mathbb{N} \setminus \{0,1\} : \gcd(n^k-1,n!-1)>1\}.
\]

The first exponent for which finiteness remains unresolved is $k=3$, for which the two smallest solutions are
\[
2283\quad\text{and}\quad138645,
\]
and both produce the same prime common divisor:
\[
\gcd(2283^3-1,2283!-1)
=
\gcd(138645^3-1,138645!-1)
=
140929.
\]
They are connected by
\[
2283+138645+1=140929,
\]
and are the two nontrivial cubic roots of unity modulo $140929$. This motivated the question in \cite{RipaCubicMO} of whether they were the only solutions, or at least whether every nontrivial gcd in the cubic case had to equal $140929$.

A stronger and more persistent feature was already present. Since $140929$ divides $2283^3-1$, it divides $2283^{3 \cdot m}-1$ for every positive integer $m$; the same is true for $138645$. Thus, each of these integers remains a solution for every exponent divisible by $3$. Extending this search from the cubic case to arbitrary odd primes led the author to the notion of a \emph{universal invariant}. More precisely, for a fixed odd prime $p$, an integer $n \geq 2$ is called a universal invariant if
\[
\gcd(n^{p \cdot m}-1,n!-1)>1
\]
holds for every integer $m \geq 1$.

Universal invariants are arithmetically interesting in their own right since their defining condition involves a persistent factorial congruence and a root of unity modulo a prime larger than the invariant itself. Moreover, the smallest universal invariants associated with the odd primes are highly irregular \cite{OEISA395171} and, in the cubic case, the two nontrivial cubic roots modulo a witness prime $q$ may or may not both satisfy the factorial congruence $n!\equiv1\pmod q$.
Universal invariants also provide a particularly simple application to the cardinality question that motivated the broader investigation: if an exponent contains many distinct odd prime factors, we may retain one universal invariant belonging to each of them; this produces an explicit set of guaranteed solutions whose size grows with the number of such prime factors.

The question of whether every odd prime admits a universal invariant was posed in~\cite{RipaUniversalMO}. In the same discussion, the author also proposed using universal invariants associated with distinct odd prime factors of the exponent to obtain lower bounds for $|S_k^{-1,-1}|$ growing with the number of such factors; the initial formulation assumed the existence of at least two distinct universal invariants for each odd prime. The cubic case was already known, and Jonathan Love settled all the remaining cases \cite{LoveMO} by proving that at least one universal invariant exists for every odd prime $p\geq5$.
Love's construction selects a prime divisor $q$ of $\frac{2^p+1}{3}$ and takes $n=q-2$, although this need not be the smallest universal invariant associated with $p$. Wilson's theorem then forces $n!\equiv1\pmod q$, while $n\equiv-2\pmod q$ forces $n^p\equiv1\pmod q$.

\begin{lemma}\label{lem:cubic-existence}
The integer $2283$ is a universal invariant for $p=3$.
\end{lemma}

\begin{proof}
A direct computation gives
\[
\gcd(2283^3-1,2283!-1)=140929>1.
\]
Since
\[
2283^3-1\mid 2283^{3 \cdot m}-1
\]
for every positive integer $m$, the divisor $140929$ divides both $2283^{3 \cdot m}-1$ and $2283!-1$ for every $m \geq 1$. Hence,
\[
\gcd(2283^{3 \cdot m}-1,2283!-1)>1
\]
holds for every positive integer $m$.
\end{proof}

In this paper, we show that one universal invariant can be chosen for each odd prime so that the resulting invariants are pairwise distinct.
If $\oddomega(k)$ denotes the number of distinct odd prime divisors of the positive integer $k$ and $\omega(k)$ the total number of distinct prime divisors, it follows that
\[
|S_k^{-1,-1}|\geq\oddomega(k)\geq \omega(k)-1.
\]
Since $\oddomega(k)$ is unbounded, so are the cardinalities $|S_k^{-1,-1}|$ as $k$ varies. This proves constructively that we can force arbitrarily many solutions while keeping both signs negative, but it remains unknown whether $S_k^{-1,-1}$ is finite or infinite for any fixed $k\geq3$.

The lower bound above is intentionally conservative; it counts only a small, rigid family of solutions provided by universal invariants. Computations for composite exponents reveal many further solutions that are not universal invariants for any individual prime divisor of the exponent, and Corollary~\ref{cor:direct-bound} later gives a stronger general lower bound by a construction not restricted to universal invariants. No attempt is made here to classify those additional solutions or to determine the exact cardinality of $S_k^{-1,-1}$. Our purpose is instead to develop the theory of the rarer universal invariants and use their persistence as a clean constructive tool.

The cubic search was subsequently extended by Aldo Roberto Pessolano using a parallel remainder-tree algorithm. His computation produced the third smallest solution
\[
27735920,
\]
with
\[
\gcd(27735920^3-1,27735920!-1)=54003043,
\]
and excluded any fourth solution below $1005000000$ \cite{PessolanoMO}.
This answered the MathOverflow question in \cite{RipaCubicMO} negatively and showed that the first pair does not exhaust the cubic phenomenon.

The algebraic component of the problem is largely classical and the existence of a common prime divisor reduces to the study of nontrivial roots of unity in the multiplicative group of a finite field, whose structure is well understood. The genuine difficulty lies instead in determining whether such roots also satisfy the factorial congruence \mbox{$n!\equiv1\pmod q$}.

The distribution of factorial residues modulo primes has been studied extensively, although it remains far from being understood. Our problem concerns factorial residues evaluated at a highly structured subset of indices, namely the nontrivial $p$-th roots of unity modulo a prime $q$.

Some of the relevant results are reviewed in Section~\ref{sec:factorial-residues}. In particular, Banks, Luca, Shparlinski, and Stichtenoth \cite{BanksLucaShparlinskiStichtenoth} established quantitative results on the residue classes omitted by factorials modulo primes; Klurman and Munsch \cite{KlurmanMunsch} later obtained an average analogue.

Although these results do not directly address the root-of-unity constraint considered here, they provide a natural mathematical context for the heuristic argument proposed by Will Sawin in the MathOverflow discussion \cite{RipaCubicMO}, which is discussed in Subsection~\ref{sec:cubic-heuristic}.

More broadly, the literature reviewed in Section~\ref{sec:factorial-residues} places the present investigation within the active study of factorial residues modulo primes.

\section{Preliminary results}

We collect here a few elementary results that will be used throughout the rest of the paper.

\begin{lemma}\label{lem:q-greater-n}
Let $n\geq2$ and $h\geq1$ be integers. Let $q$ be a prime divisor of
\[
\gcd(n^h-1,n!-1).
\]
Then, $q>n$.
\end{lemma}

\begin{proof}
If $q\leq n$, then $q\mid n!$, contradicting $q\mid n!-1$.
\end{proof}

The following general prime-witness reduction was pointed out in a comment to the MathOverflow discussion \cite{RipaCardinalityMO}.

\begin{proposition}\label{prop:general-witness}
Let $n\geq2$ and $k\geq1$ be integers, let $\lambda_1,\lambda_2\in\{-1,+1\}$, and let $q=n+r$ be prime, with $r\geq1$. Then,
\[
q\mid\gcd(n^k+\lambda_1,n!+\lambda_2)
\]
if and only if
\[
(-r)^k+\lambda_1\equiv0\pmod q
\]
and
\[
(r-1)!\equiv\lambda_2(-1)^{r-1}\pmod q.
\]
\end{proposition}

\begin{proof}
Since $n\equiv-r\pmod q$, the first equivalence is immediate. Moreover, with the empty product understood as $1$ when $r=1$,
\[
(q-1)!
=
n!\prod_{j=1}^{r-1}(n+j)
\equiv
n!\cdot(-1)^{r-1}\cdot(r-1)!
\pmod q.
\]
Wilson's theorem then gives the second equivalence.
\end{proof}

\begin{lemma}\label{lem:roots}
Let $h\geq1$ and let $q$ be prime. The equation $x^h=1$ has exactly $\gcd(h,q-1)$ solutions in $\F_q^\ast$.
\end{lemma}

\begin{proof}
Since $\F_q^\ast$ is cyclic of order $q-1$ (see, e.g., \cite[pp.~239--240]{Aluffi}), write $x=g^u$, where $g$ is a generator. Then, $x^h=1$ if and only if $q-1\mid h \cdot u$. Thus, there are exactly $\gcd(h,q-1)$ residue classes $u\pmod{q-1}$ satisfying this condition.
\end{proof}

\begin{corollary}\label{cor:p-roots}
Let $p$ and $q$ be primes. If $p\mid q-1$, then $x^p=1$ has exactly $p$ solutions in $\F_q^\ast$, while if $p\nmid q-1$, it has only the trivial solution $x=1$.
\end{corollary}

\begin{remark}
When $q=p$, we have $x^p-1=(x-1)^p$ in $\F_p[x]$. Thus, there is still only one distinct root, namely $1$, although it has multiplicity $p$.
\end{remark}

\begin{proposition}\label{prop:common-prime}
Let $p$ be an odd prime, $n\geq2$, and suppose that a prime $q$ divides
\[
\gcd(n^p-1,n!-1).
\]
Then, $q>n$ and $n\not\equiv1\pmod q$. Moreover, $\operatorname{ord}_q(n)=p$ and $q\equiv1\pmod{2 \cdot p}$.
\end{proposition}

\begin{proof}
By Lemma~\ref{lem:q-greater-n}, we have $q>n$. Hence, $1<n<q$, so the residue class of $n$ belongs to $\F_q^\ast\setminus\{1\}$. Since $n^p\equiv1\pmod q$, the multiplicative order $\operatorname{ord}_q(n)$ divides $p$. As $p$ is prime and $n\not\equiv1\pmod q$, we must have $\operatorname{ord}_q(n)=p$.

Consequently, $p\mid q-1$. Since $q>n\geq2$, the prime $q$ is odd, so $2\mid q-1$. As $p$ is odd, it follows that $2 \cdot p\mid q-1$. Therefore, $q\equiv1\pmod{2 \cdot p}$.
\end{proof}

\section{Universal invariants}

\begin{definition}\label{def:universal-invariant}
Let $p$ be an odd prime. An integer $n\geq2$ is called a ``universal invariant'' for $p$ if
\[
\gcd(n^{p \cdot m}-1,n!-1)>1
\]
holds for every positive integer $m$.
\end{definition}

\begin{proposition}\label{prop:k-one}
Let $p$ be an odd prime and let $n\geq2$ be an integer. The following are equivalent:
\begin{enumerate}[label=\textup{(\roman*)}]
\item $\gcd(n^{p \cdot m}-1,n!-1)>1$ for every integer $m\geq1$;
\item $\gcd(n^p-1,n!-1)>1$;
\item there exists a prime $q>n$ such that $\operatorname{ord}_q(n)=p$ and $n!\equiv1\pmod q$.
\end{enumerate}
\end{proposition}

\begin{proof}
Condition (i) implies (ii) by taking $m=1$. Conversely,
\[
n^{p \cdot m}-1=(n^p)^m-1
=(n^p-1)\cdot\left(1+n^p+\cdots+n^{p \cdot (m-1)}\right).
\]
Thus, $n^p-1\mid n^{p \cdot m}-1$ for every $m \geq 1$, so every common divisor of $n^p-1$ and $n!-1$ is also a common divisor of both $n^{p \cdot m}-1$ and $n!-1$.

Now, suppose that (ii) holds. Let $q$ be any prime divisor of $\gcd(n^p-1,n!-1)$. By Proposition~\ref{prop:common-prime}, we have $q>n$ and $\operatorname{ord}_q(n)=p$. Moreover, $q\mid n!-1$, and hence $n!\equiv1\pmod q$. Thus, (iii) holds.

Finally, if (iii) holds, then $\operatorname{ord}_q(n)=p$ implies $n^p\equiv1\pmod q$. Together with $n!\equiv1\pmod q$, this gives $q\mid\gcd(n^p-1,n!-1)$, so (ii) holds.
\end{proof}

\begin{theorem}[Mantovanelli \cite{MantovanelliMO}]\label{thm:p-two}
For every integer $n\geq2$,
\[
\gcd(n^2-1,n!-1)=1.
\]
\end{theorem}

\begin{proof}
Suppose that a prime $q$ divides both $n^2-1$ and $n!-1$. By Lemma~\ref{lem:q-greater-n}, $q>n$. Since $q\mid(n-1)\cdot(n+1)$ and $0<n-1<q$, we have $q\nmid n-1$. Hence, $q\mid n+1$. Since $0<n+1<2 \cdot q$, it follows that $n+1=q$ (so that $n=q-1$). Finally, Wilson's theorem gives $n!=(q-1)!\equiv-1\pmod q$, contradicting $q\mid n!-1$.
\end{proof}

\begin{remark}
No universal invariant exists for $p=2$, in the sense of Definition~\ref{def:universal-invariant} with $p=2$ allowed (equivalently, $|S_2^{-1,-1}|=0$).
\end{remark}

\subsection{Explicit construction and a cardinality consequence}

\begin{theorem}[Love \cite{LoveMO}]\label{thm:explicit-construction}
Let $p\geq5$ be an odd prime and let $q$ be any prime divisor of $\frac{2^p+1}{3}$. Then, $n=q-2$ is a universal invariant for $p$. More precisely,
\[
q\mid\gcd(n^{p \cdot m}-1,n!-1)
\]
for every positive integer $m$.
\end{theorem}

\begin{proof}
Since $p\geq5$ is prime, we have $p\equiv\pm1\pmod6$. Hence, the integer $2^p+1$ is odd, divisible by $3$, but not by $3^2$.

Thus, every prime divisor $q$ of $\frac{2^p+1}{3}$ satisfies $q\geq5$.
Let $n:=q-2$. Since $p$ is odd and $q\mid2^p+1$, Proposition~\ref{prop:general-witness}, applied with $k=p$, $r=2$, and $\lambda_1=\lambda_2=-1$, yields
\[
q\mid\gcd(n^p-1,n!-1).
\]
Proposition~\ref{prop:k-one} now completes the proof.
\end{proof}

\begin{corollary}\label{cor:existence}
A universal invariant exists for every odd prime.
\end{corollary}

\begin{proof}
For $p=3$, the integer $2283$ is a universal invariant. For every odd prime $p\geq5$, apply Theorem~\ref{thm:explicit-construction}.
\end{proof}

In view of Corollary~\ref{cor:existence}, the following minimum is well defined for every odd prime.

\begin{definition}
For any odd prime $p$, let
\[
\inv_p := \min\{\,n \geq 2:\gcd(n^p-1,n!-1)>1\,\}.
\]
We call $\inv_p$ the ``minimal universal invariant'' for $p$.
\end{definition}

\begin{corollary}\label{cor:upper-bound}
Let $p\geq5$ be prime and let $q$ be the least prime divisor of $\frac{2^p+1}{3}$. Then,
\[
\inv_p\leq q-2\leq\frac{2^p+1}{3}-2.
\]
\end{corollary}

\begin{lemma}\label{lem:distinct-q}
Let $p$ and $s$ be distinct odd primes. No prime $q\geq5$ can divide both $2^p+1$ and $2^s+1$.
\end{lemma}

\begin{proof}
If $q$ divides $2^p+1$, then the multiplicative order of $2$ modulo $q$ divides $2 \cdot p$ but not $p$. Since $p$ is prime, the order is therefore either $2$ or $2 \cdot p$. The former is impossible for $q\geq5$ since it would imply $q\mid 2^2-1=3$. Hence, the order is $2 \cdot p$. Similarly, if $q$ divides $2^s+1$, then the order is $2 \cdot s$. Thus, $2 \cdot p=2 \cdot s$, and this contradicts $p\neq s$.
\end{proof}

\begin{theorem}\label{thm:omega-bound}
For every integer $k\geq2$,
\[
|S_k^{-1,-1}|\geq\oddomega(k).
\]
Consequently, the cardinalities $|S_k^{-1,-1}|$ are unbounded as $k$ varies.
\end{theorem}

\begin{proof}
Let $t:=\oddomega(k)$, and let $p_1,\ldots,p_t$ be the distinct odd prime divisors of $k$. For each odd prime divisor $p_i$ of $k$, define $n_i$ as follows: if $p_i=3$, let $n_i=2283$. Otherwise, let $q_i$ be any prime divisor of $\frac{2^{p_i}+1}{3}$ and put $n_i=q_i-2$. Lemma~\ref{lem:cubic-existence} shows that $2283$ is a universal invariant for $3$ while Theorem~\ref{thm:explicit-construction} gives a universal invariant for every odd prime $p_i \geq 5$. Since $p_i\mid k$, we have $n_i\in S_k^{-1,-1}$.
Lemma~\ref{lem:distinct-q} shows that the $n_i$ corresponding to $p_i\geq5$ are pairwise distinct. Furthermore, none of them equals $2283$ (since this would require the prime $q_i$ to equal $2285=5 \cdot 457$). Hence, $|S_k^{-1,-1}|\geq t=\oddomega(k)$.
\end{proof}

\begin{remark}
Theorem~\ref{thm:omega-bound} is obtained by constructing one explicit universal invariant for each distinct odd prime divisor of $k$. The random-residue heuristic discussed in Subsection~\ref{sec:cubic-heuristic} extends naturally to every fixed odd prime $p$ and suggests that the expected number of primes admitting such a root diverges. Although the underlying independence assumptions are unlikely to hold exactly, the heuristic still supports the possibility that infinitely many universal invariants exist for every fixed odd prime.
\end{remark}

\subsection{Minimal universal invariants}

The values of the minimal universal invariants for odd primes up to $59$ currently recorded in OEIS sequence A395171 \cite{OEISA395171} are listed in Table~\ref{tab:minimal}.

\begin{table}[H]
\centering
\begin{tabular}{|c|c|}
\hline
$p$ & $\inv_p$\\
\hline\hline
3 & 2283\\
\hline
5 & 9\\
\hline
7 & 41\\
\hline
11 & 4\\
\hline
13 & 15\\
\hline
17 & 43689\\
\hline
19 & 174761\\
\hline
23 & 236\\
\hline
29 & 29\\
\hline
31 & 357913941\\
\hline
37 & 25\\
\hline
41 & 41\\
\hline
43 & 64\\
\hline
47 & 141\\
\hline
53 & 53\\
\hline
59 & 2831\\
\hline
\end{tabular}
\caption{Minimal universal invariants for odd primes up to $59$.}
\label{tab:minimal}
\end{table}

Love's construction can be sharp, but it may also give a much larger value. For $p=11$, $\frac{2^{11}+1}{3}=683$ is prime, so the construction gives $681$; nevertheless, $\inv_{11}=4$. Indeed,
\[
4^{11}\equiv1\pmod{23}
\qquad\text{and}\qquad
4!-1=23.
\]
At the opposite extreme, $\inv_{31}=357913941$ (its minimality was confirmed by Pessolano through an exhaustive computation \cite{Pessolano31}).

\subsection{Persistent common divisors}

\begin{definition}
Let $p$ be an odd prime and let $n$ be a universal invariant for $p$. Define
\[
D_p(n):=\gcd\left\{\,\gcd(n^{p \cdot m}-1,n!-1):m\geq1\,\right\},
\]
where the gcd of the infinite family denotes the greatest positive integer dividing every member of the family. This is well defined since every member divides the fixed integer $n!-1$.
\end{definition}

\begin{proposition}\label{prop:persistent-divisor}
For every universal invariant $n$ for $p$,
\[
D_p(n)=\gcd(n^p-1,n!-1)>1.
\]
\end{proposition}

\begin{proof}
Since $n^p-1\mid n^{p \cdot m}-1$ for every $m \geq 1$, we have
\[
\gcd(n^p-1,n!-1)\mid\gcd(n^{p \cdot m}-1,n!-1)
\]
for every $m \geq 1$. Hence,
\[
\gcd(n^p-1,n!-1)\mid D_p(n).
\]
Conversely, the family defining $D_p(n)$ contains the case $m=1$, so
\[
D_p(n)\mid\gcd(n^p-1,n!-1).
\]
Therefore, the two quantities are equal.
\end{proof}

For consecutive odd primes $p$, the values of $D_p(\inv_p)$ begin
\[
140929,121,43,23,53,43691,174763,277,59,715827883,149,83,431,283,107,2833,\ldots
\]
(they are recorded as OEIS sequence A393224 \cite{OEISA393224}).

It is worth noting that the persistent common divisor need not be prime or even squarefree (e.g., $D_5(\inv_5)=121=11^2$ so that $11^2\mid\gcd(9^{5 \cdot m}-1,9!-1)$ for every positive integer $m$).

The statement above does not imply that the entire gcd is independent of $m$: additional prime factors may occur for particular multipliers. More generally, if $\ell$ is a prime divisor of $n!-1$, then necessarily $\ell>n$, and hence
\[
\ell\mid n^{p \cdot m}-1
\quad\Longleftrightarrow\quad
\operatorname{ord}_{\ell}(n)\mid p \cdot m.
\]
Thus, the occurrence of a prime factor $\ell$ in the gcd at multiplier $m$ is governed exactly by this order-divisibility condition. For example, direct computation for each integer $1 \leq m \leq 2167$ gives
\[
\gcd(2283^{3 \cdot m}-1,2283!-1)=140929,
\]
whereas
\[
\gcd(2283^{3 \cdot 2168}-1,2283!-1)
=611209073=4337 \cdot 140929.
\]
Indeed,
\[
\operatorname{ord}_{4337}(2283)=2168,
\]
and, since $\gcd(2168,3)=1$, the factor $4337$ occurs precisely when $2168 \mid m$, explaining its first appearance at $m=2168$.

\section{The cubic case}

\begin{proposition}\label{prop:cubic-reduction}
For every integer $n\geq2$,
\[
\gcd(n^3-1,n!-1)>1
\]
if and only if
\[
\gcd(n^2+n+1,n!-1)>1.
\]
\end{proposition}

\begin{proof}
Since
\[
n^3-1=(n-1)\cdot(n^2+n+1),
\]
it suffices to show that
\[
\gcd(n-1,n!-1)=1.
\]
Otherwise, a prime divisor of $\gcd(n-1,n!-1)$ would be at most $n-1$, whereas Lemma~\ref{lem:q-greater-n}, applied with $h=1$, implies that it must be greater than $n$.
\end{proof}

\begin{definition}\label{paired}
Let $q\equiv1\pmod3$ be prime and let $a,b\in\{1,\ldots,q-1\}$ be the two nontrivial cubic roots of unity modulo $q$. We say that $a$ and $b$ are $q$-paired if
\[
a!\equiv b!\equiv1\pmod q.
\]
\end{definition}

The parity phenomenon underlying the following proposition was observed by Pessolano in \cite{PessolanoMO}: we record it here in a general form.

\begin{proposition}\label{prop:parity-criterion}
Let $q\equiv1\pmod3$ be prime and let $a,b\in\{1,\ldots,q-1\}$ be the two nontrivial cubic roots of unity in $\F_q$. Then,
\[
a+b+1=q
\]
and
\[
a! \cdot b!\equiv(-1)^{a+1}\pmod q.
\]
In particular, if $a!\equiv1\pmod q$, then
\[
b!\equiv(-1)^{a+1}\pmod q.
\]
Consequently, under the hypothesis $a!\equiv1\pmod q$,
\[
q\mid\gcd(b^3-1,b!-1)
\]
if and only if $a$ is odd.
\end{proposition}

\begin{proof}
The roots $a$ and $b$ satisfy
\[
x^2+x+1=0,
\]
so
\[
a+b\equiv-1\pmod q.
\]
Since
\[
2 \leq a+b \leq 2 \cdot q-2,
\]
it follows that
\[
a+b=q-1.
\]
Hence,
\[
b=q-a-1.
\]
Now,
\[
(q-1)!
=
b!\prod_{j=b+1}^{q-1}j
=
b!\prod_{j=1}^{a}(q-j)
\equiv
b! \cdot (-1)^a \cdot a!
\pmod q.
\]
By Wilson's theorem,
\[
-1\equiv(q-1)!\equiv b! \cdot (-1)^a \cdot a!\pmod q.
\]
Therefore,
\[
a! \cdot b!\equiv(-1)^{a+1}\pmod q.
\]
If
\[
a!\equiv1\pmod q,
\]
then
\[
b!\equiv(-1)^{a+1}\pmod q.
\]
Since $q$ is odd,
\[
(-1)^{a+1}\equiv1\pmod q
\]
if and only if $a$ is odd. Moreover, since $b$ is a cubic root of unity modulo $q$, we have
\[
q\mid b^3-1.
\]
Hence,
\[
q\mid\gcd(b^3-1,b!-1)
\]
if and only if $a$ is odd.
\end{proof}

The first two cubic universal invariants, $2283$ and $138645$, are the two nontrivial cubic roots modulo $140929$ and are $140929$-paired (see Definition~\ref{paired}):
\[
2283+138645+1=140929.
\]
The third smallest cubic universal invariant, found by Pessolano \cite{PessolanoMO}, is
\[
27735920,
\]
with
\[
\gcd(27735920^3-1,27735920!-1)=54003043.
\]
The other nontrivial cubic root modulo $54003043$ is $26267122$.

Since $27735920$ is even, Proposition~\ref{prop:parity-criterion} gives
\[
26267122!\equiv-1\pmod{54003043},
\]
and consequently the two nontrivial cubic roots modulo $54003043$ are not $54003043$-paired.

Finally, the exhaustive search described in~\cite{PessolanoMO} shows that no fourth cubic solution occurs below $1005000000$.

\subsection{Heuristic for the cubic case}\label{sec:cubic-heuristic}

For every prime $q\equiv1\pmod3$, there are exactly two nontrivial cubic roots modulo $q$. Treating the factorial residue associated with each nontrivial cubic root as uniformly distributed in $\F_q^\ast$, each root has heuristic probability $\frac{1}{q-1}$ of satisfying $n!\equiv1\pmod q$. The two events should not, however, be regarded as independent, in view of Proposition~\ref{prop:parity-criterion}.
In a comment to \cite{RipaCubicMO}, Will Sawin observed that this simple random heuristic suggests an expected number of contributing primes up to $X$ of order $\log\log X$.
This suggests the following conjecture.

\begin{conjecture}
There exist infinitely many positive integers $n$ such that
\[
\gcd(n^3-1,n!-1)>1.
\]
\end{conjecture}

\begin{remark}
The heuristic above should be regarded only as a first approximation. It treats the relevant factorial residues as essentially random and does not take into account possible arithmetic restrictions that may cause different residue classes to contribute with different frequencies.
\end{remark}

\subsection{Factorial residues modulo primes}\label{sec:factorial-residues}

Let
\[
F_q:=\{\,n!\bmod q:\ 1 \leq n<q\,\}\subseteq\F_q^\ast
\]
denote the set of factorial residues modulo the prime $q$.

Thus, $q-|F_q|$ counts the residue classes in $\F_q$ omitted by the factorial sequence, including the always-missing class $0$.

Our problem asks whether the residue $1$ occurs at a point which is simultaneously a nontrivial $p$-th root of unity. It therefore belongs to the broader study of the distribution of factorial residues modulo primes.

Cobeli, V\^aj\^aitu, and Zaharescu studied factorial residues modulo primes and introduced a probabilistic model in which the average proportion of missing values tends to $\frac{1}{e}$ \cite{CobeliVajaituZaharescu}. Subsequent work showed that the value set $F_q$ possesses significant arithmetic structure. Banks, Luca, Shparlinski, and Stichtenoth proved that there exist infinitely many primes $q$ such that
\[
q-|F_q|
\gg
\frac{\log\log q}{\log\log\log q},
\]
so the number of residue classes omitted by $1!,\ldots,(q-1)!$ is unbounded \cite{BanksLucaShparlinskiStichtenoth}.
Klurman and Munsch later established an average analogue by proving that
\[
\frac1{\pi(x)}
\sum_{q \leq x}\bigl(q-|F_q|\bigr)
\gg
\frac{\log\log x}{\log\log\log x},
\]
showing that the average number of omitted residue classes also tends to infinity \cite{KlurmanMunsch}.

In 2024, Grebennikov, Sagdeev, Semchankau, and Vasilevskii proved the complementary lower bound
\[
|F_q|\geq (\sqrt{2}+o(1))\sqrt q
\]
as $q\to\infty$ through the primes \cite[Corollary~2.1]{GrebennikovSagdeevSemchankauVasilevskii}.

More recently, Hu broke the square-root barrier by proving
\[
|F_q|\gg q^{8/15}
\]
for every prime $q$ \cite{Hu}.

Further aspects of factorial residues have been investigated through random models for the missing residue classes \cite{BroughanBarnett,CobeliPitropZaharescu}, character sums, exponential sums, and congruences involving factorials and products of factorials \cite{GaraevLuca,GaraevLucaShparlinski,GaraevLucaShparlinskiExp}.

Finally, since $1!\equiv(q-2)!\equiv1\pmod q$ for every prime $q\geq5$, we have $q-|F_q|\geq2$. Equality holds precisely when this forced repetition is the only one, or equivalently when $2!,3!,\ldots,(q-1)!$ are pairwise distinct modulo $q$. Primes $q>5$ with this property are usually called socialist primes. Trudgian excluded such primes below $10^9$ \cite{Trudgian}; Andreji\'c and Tatarevi\'c subsequently extended the computation to $10^{11}$ \cite{AndrejicTatarevic}, and Andreji\'c, Bostan, and Tatarevi\'c further extended it to $2^{40}$ \cite{AndrejicBostanTatarevic}.
More recently, a preprint by Abramov claims a complete negative solution to the problem \cite{Abramov}.

Despite these advances, none of the existing results directly addresses the present problem: here, the residue is constrained to lie among the nontrivial $p$-th roots of unity, so possible correlations between the root condition $n^p\equiv1\pmod q$ and the factorial congruence $n!\equiv1\pmod q$ remain essentially unexplored.

\subsection{Comparison with the original four-sign family}

The behaviour of the generalized family
\[
\gcd(n^k+\lambda_1,n!+\lambda_2)>1,
\qquad \lambda_1,\lambda_2\in\{-1,+1\},
\]
strongly depends on the choice of signs and on the parity of $k$.

\begin{proposition}
Let $q$ be an odd prime and set $n:=q-1$. Then,
\[
\gcd(n^k-1,n!+1)>1
\]
holds for every even positive integer $k$, while
\[
\gcd(n^k+1,n!+1)>1
\]
holds for every odd positive integer $k$.
\end{proposition}

\begin{proof}
Apply Proposition~\ref{prop:general-witness} with $r=1$.
\end{proof}

Thus, $S_k^{-1,+1}$ is infinite for every even $k$, while $S_k^{+1,+1}$ is infinite for every odd $k$; in fact, every odd prime occurs as a common divisor in the corresponding case.

A complementary construction applies when $\lambda_2=-1$.

\begin{proposition}
Let $q\geq5$ be prime and set $n:=q-2$. Then,
\[
q\mid\gcd(n^k-1,n!-1)
\]
if and only if
\[
q\mid 2^k-(-1)^k,
\]
whereas
\[
q\mid\gcd(n^k+1,n!-1)
\]
if and only if
\[
q\mid 2^k+(-1)^k.
\]
\end{proposition}

\begin{proof}
Apply Proposition~\ref{prop:general-witness} with $r=2$ and $\lambda_2=-1$.
\end{proof}

Thus, the construction $n=q-2$ yields solutions in $S_k^{-1,-1}$ whenever $q\geq5$ is a prime divisor of $2^k-(-1)^k$, and in $S_k^{+1,-1}$ whenever $q\geq5$ is a prime divisor of $2^k+(-1)^k$.

Moreover, for every $k\geq2$, the integer $2^k+(-1)^k$ is odd, greater than $1$, and not divisible by $3$. Hence, it has a prime divisor $q\geq5$, so $S_k^{+1,-1}\neq\varnothing$ for every $k\geq2$.

Let $\omega_{\geq5}(m)$ denote the number of distinct prime divisors of $m$ that are at least $5$.

\begin{corollary}\label{cor:direct-bound}
For every integer $k\geq2$,
\[
|S_k^{-1,-1}|
\geq
\omega_{\geq5}\bigl(2^k-(-1)^k\bigr)
+\mathbf{1}_{3\mid k},
\]
where $\mathbf{1}_{3\mid k}=1$ if $3\mid k$ and $0$ otherwise.
\end{corollary}

\begin{proof}
Every distinct prime $q\geq5$ dividing $2^k-(-1)^k$ yields the distinct solution $n=q-2$ by the preceding proposition. If $3\mid k$, then $2283\in S_k^{-1,-1}$ by Lemma~\ref{lem:cubic-existence}. This solution is distinct from all those of the form $q-2$, since $2283+2=2285=5 \cdot 457$ is not prime.
\end{proof}

\begin{remark}
An application of the classical Bang--Zsigmondy theorem \cite{Bang,Zsigmondy} yields the further bound
\[
|S_k^{-1,-1}|
\geq
\begin{cases}
\tau(k)-1 & \text{if $k$ is odd}\\
\tau(k)-2 & \text{if $k$ is even}
\end{cases},
\]
where $\tau(k)$ denotes the number of positive divisors of $k$.

Indeed, if $k$ is odd, then for every divisor $d>1$ of $k$, the Bang--Zsigmondy theorem applied to $2^{2 \cdot d}-1$ gives, except when $d=3$, a prime $q_d$ with $\operatorname{ord}_{q_d}(2)=2 \cdot d$. Hence, $q_d\mid 2^d+1$, and since $k/d$ is odd, also $q_d\mid 2^k+1$; the primes obtained for distinct $d$ are distinct. When $3\mid k$, the exceptional divisor $d=3$ is compensated by the solution $2283$.

If $k$ is even, applying the theorem to $2^d-1$ for divisors $d>1$ of $k$ gives distinct prime divisors of $2^k-1$. The case $d=2$ yields only the prime $3$, which is not admissible in the $q=n+2$ construction, while the exceptional case $d=6$, when present, is again compensated by $2283$.

We do not pursue this refinement further here, since our main interest is in the structure and persistence of universal invariants rather than in optimizing general lower bounds for $|S_k^{-1,-1}|$.
\end{remark}

The main results of this subsection are summarized in Table~\ref{tab:four-sign-family}.

\begin{table}[ht]
\centering
\renewcommand{\arraystretch}{1.25}
\begin{tabular}{|>{\centering\arraybackslash}m{0.28\textwidth}|
                >{\centering\arraybackslash}m{0.20\textwidth}|
                >{\centering\arraybackslash}m{0.40\textwidth}|}
\hline
CONDITION & RESULT & MAIN ARGUMENT\\
\hline\hline

$\gcd(n^2-1,n!-1)>1$ &
No solutions &
Mantovanelli's theorem.\\
\hline

$\gcd(n^k-1,n!-1)>1$, $k\geq3$ &
Cardinalities unbounded as $k$ varies &
Theorem~\ref{thm:omega-bound} gives $|S_k^{-1,-1}|\geq\oddomega(k)$ using universal invariants; Corollary~\ref{cor:direct-bound} gives a stronger general lower bound.\\
\hline

$\gcd(n^k-1,n!+1)>1$, $k\geq2$ even &
Infinitely many solutions &
Take $n=q-1$ for any odd prime $q$.\\
\hline

$\gcd(n^k+1,n!-1)>1$, $k\geq2$ &
Nonempty &
The integer $2^k+(-1)^k$ has a prime divisor $q\geq5$; take $n=q-2$.\\
\hline

$\gcd(n^2+1,n!+1)>1$ &
Open problem &
The original MathOverflow problem.\\
\hline

$\gcd(n^k+1,n!+1)>1$, $k\geq3$ odd &
Infinitely many solutions &
Take $n=q-1$ for any odd prime $q$.\\
\hline

\end{tabular}
\caption{Representative behaviours in the four-sign family $\gcd(n^k+\lambda_1,n!+\lambda_2)>1$ for integer $k\geq2$.}\label{tab:four-sign-family}
\end{table}

\section{Open problems}

This paper focuses on the negative-sign family $\gcd(n^k-1,n!-1)>1$, where the notion of universal invariants naturally arises.

Despite the results established in the preceding sections, several interesting questions remain open.

\begin{problem}\label{problem:fixed-k}
For each fixed integer $k\ge3$, is $S_k^{-1,-1}$ finite or infinite?
\end{problem}

The following problem is the original four-sign formulation that initiated the present investigation \cite{RipaCardinalityMO}.

\begin{problem}
For every positive integer $c$, do there exist an integer $k\geq2$ and signs $\lambda_1,\lambda_2\in\{-1,+1\}$ such that
\[
c<|S_k^{\lambda_1,\lambda_2}|<\infty?
\]
\end{problem}

By Proposition~\ref{prop:k-one}, the following problem is the restriction of Open Problem~\ref{problem:fixed-k} to odd prime exponents.

\begin{problem}
For a fixed odd prime $p$, are there infinitely many universal invariants for $p$?
\end{problem}

For $p=3$, the heuristic discussed in Subsection~\ref{sec:cubic-heuristic} is consistent with an affirmative answer.

\begin{problem}
Can substantially sharper bounds for $\inv_p$ be obtained?
\end{problem}

\begin{problem}
Can the distribution of factorial residues over the nontrivial $p$-th roots of unity be described more precisely than by a uniform random model?
\end{problem}

Finally, we recall the question that initiated this line of investigation and remains open (see \cite{SmallerMO,MounierMO}): does there exist an integer $n\geq2$ such that $\gcd(n^2+1,n!+1)>1$?

\section*{Acknowledgements}

The author thanks Alex Casarotti for his continuous support and several useful discussions concerning the conjectured infinitude of universal invariants. He is grateful to Francesco Santilli for an independent review of an earlier version of the manuscript; to Marco Mantovanelli for proving the quadratic negative-sign case; to Jonathan Love for his construction of universal invariants for odd primes $p\geq5$; to Will Sawin for the heuristic argument in the cubic case; and to Aldo Roberto Pessolano for his extensive computations in the cubic case, including the discovery of the third smallest cubic universal invariant and the exhaustive verification that $\inv_{31}=357913941$. He also thanks the other MathOverflow contributors whose questions, answers, and comments helped shape the present work.

\section*{Declaration on the use of AI}

During the preparation of this manuscript, the author used ChatGPT (OpenAI, GPT-5.5 and GPT-5.6) as an interactive writing, editing, and analytical assistant. The system was employed to discuss mathematical ideas, review the presentation of proofs, improve the organization of the manuscript, suggest alternative expository formulations, and assist in the preparation of the \LaTeX{} source.

All mathematical results presented in this article originate either from the cited sources or from the author's own work. The author independently verified the mathematical arguments, references, and final formulations presented here; computational results due to external contributors are reported with attribution to the corresponding sources. The author assumes full responsibility for the contents of this manuscript.

\end{document}